\documentclass[oneside,reqno]{amsart}
\usepackage[greek,english]{babel}
\usepackage{amsthm}
\usepackage{amsbsy}
\usepackage{amsfonts}
\usepackage{graphicx}
\newtheorem{pro}{Proposition}
\newtheorem{rem}{Remark}

\begin{document}

\title[Counterexamples to Modica's gradient estimate for systems]{Counterexamples to Modica's gradient estimate for systems arising in multi-phase transitions}

\author{Christos Sourdis}
\address{Department of Mathematics,  University of
Turin,   via Carlo Alberto 10,
20123, Turin, Italy.}
\email{christos.sourdis@unito.it}
\maketitle

\begin{abstract}
We construct a plethora of one-dimensional periodic solutions for a class of semilinear elliptic systems of phase transition type which violate   Modica's gradient estimate. This complements a recent technical counterexample in \cite{smyrnelis} and is partly motivated by a recent open problem in \cite{fazly}.
\end{abstract}
 \maketitle

Consider the semilinear elliptic system
\begin{equation}\label{eqEq}
\Delta u= \nabla_u W(u),\ \ u:\mathbb{R}^n \to \mathbb{R}^m,\ \ n,m\geq 1,
\end{equation}
with $W$ sufficiently smooth and \emph{nonnegative}.

If $m=1$, it was  shown by Modica in \cite{modica} (see also \cite{cafa}) that every entire (i.e. defined in the whole space), bounded solution satisfies the pointwise gradient estimate
\begin{equation}\label{eqGB}
\frac{1}{2}|\nabla u|^2\leq W(u), \ \ x\in \mathbb{R}^n.
\end{equation}

If $m\geq 2$, the analog of the above property does not hold in general. Indeed, if it was satisfied by some entire solution (not necessarily bounded),
then the following monotonicity formula would hold:
\begin{equation}\label{eqMonotStrong}
\frac{d}{dr}\left(\frac{1}{r^{n-1}}\int_{|x|<r}^{}\left\{\frac{1}{2}|\nabla
u|^2+ W\left(u\right) \right\}dx\right)\geq 0,\ \ r>0,
\end{equation}
(see for instance \cite{alikakosBasicFacts}).
However, in the case of the Ginzburg-Landau system arising in superconductivity, where \begin{equation}\label{eqGL}
W(u)=\frac{1}{4}\left(1-|u|^2\right)^2,\end{equation}
there are bounded, entire solutions to (\ref{eqEq}) which violate the aforementioned monotonicity formula (see \cite{brezisMerle,pisante}). These solutions satisfy
\[
\frac{1}{\ln r}\int_{|x|<r}^{}\left\{\frac{1}{2}|\nabla
u|^2+ W\left(u\right) \right\}dx\to c_2,\ \ \textrm{as}\ r\to \infty, \ \ \textrm{if}\ n=m=2,
\]
and
\[
\frac{1}{ r^{n-2}}\int_{|x|<r}^{}\left\{\frac{1}{2}|\nabla
u|^2+ W\left(u\right) \right\}dx\to c_n,\ \ \textrm{as}\ r\to \infty, \ \ \textrm{if}\ n=m\geq 3,
\]
for some $c_i>0$, $i\geq 2$. In this regard, let us point out that a weak version of the monotonicity formula (\ref{eqMonotStrong}), where the exponent $n-1$
is replaced by $n-2$, holds for any solution of (\ref{eqEq}) with any $W\geq 0$ smooth (for $n\geq 2$, see for example \cite{alikakosBasicFacts,smets}). In fact, as was explained in \cite{smets}, the exponent $n-2$ is the natural one for the Ginzburg-Landau system.
Actually, in the case of the Ginzburg-Landau system with $n=1$ and $m=2$, the gradient estimate (\ref{eqGB}) is violated by the periodic solutions
\begin{equation}\label{eqper}
u_\theta(x)=\sqrt{1-\theta^2}e^{i\theta x},
\end{equation}
provided that
\begin{equation}\label{eqtheta}
\theta\in \left(0,\sqrt{\frac{2}{3}}\right).
\end{equation}
Indeed, an easy calculation (see also \cite{smyrnelis}) shows that
\begin{equation}\label{eqIneq}
\frac{1}{2}\left|u_\theta'(x) \right|^2-\frac{1}{4}\left(1-\left|u_\theta(x) \right|^2\right)^2=\frac{\theta^2}{4}(2-3\theta^2).
\end{equation}
In passing, we note that this family of periodic solutions seems to have been first observed in \cite{brezisMerle}.

From now on, we will consider system (\ref{eqEq}) with $W\geq 0$ having a finite number of global minima, which typically are assumed to be nondegenerate, appears mainly in the study of multi-phase transitions (see \cite{alikakosBasicFacts} and the references therein). But also in other contexts such as  the study of two-component Bose-Einstein condensates (see \cite{alamaARMA15} and the references therein). In that case, as explained in \cite{alikakosBasicFacts}, the number $n-1$ is the natural exponent in the denominator in (\ref{eqMonotStrong}). In fact, several properties related to (\ref{eqMonotStrong}) (with exponent $n-1$) have recently been shown to hold in \cite{AlikakosDensity} for this class of systems in the case of energy minimal solutions (in the sense of Morse). In this regard, it is of interest to know whether the analog of Modica's gradient estimate (\ref{eqGB}) holds for bounded, entire solutions to this class of systems (see the related comments in \cite{alikakosBasicFacts} and Open Problem 1 in \cite{fazly}).
It is worth noting  that a class of such systems which satisfy this property
has been provided in \cite{GhosubPass}. On the other hand, a counterexample to this property, for such systems, was provided very recently by \cite{smyrnelis}. The main idea was to start with a very particular (but smooth) periodic  function $u:\mathbb{R}\to \mathbb{R}^2$ and then build around it an axially symmetric
 $W:\mathbb{R}^2 \to \mathbb{R}$, with just two global minima, such that $u$ solves (\ref{eqEq}) and passes through the global minimizers
of $W$ with nonzero velocity. The gradient estimate (\ref{eqGB}) is clearly violated there (as well as the corresponding Liouville type theorem in \cite{cafa,modica}). In our opinion, the whole construction was rather artificial. Roughly speaking, the  curve representation of the aforementioned periodic solution on the plane consisted of two parallel line segments which were joined together smoothly with two semicircles.

Our purpose in this short note is to provide  considerably simpler counterexamples by considering the axially symmetric $W:\mathbb{R}^2 \to \mathbb{R}$ with four global minima (one on each half-axis), considered in \cite{alamaARMA15}, in the parameter regime where $W$   is  close to the Ginzburg-Landau potential (\ref{eqGL}). As the reader may have already guessed, the counterexamples will be one-dimensional periodic solutions which are smooth perturbations of the family (\ref{eqper})
with $\theta$ as in (\ref{eqtheta}) but sufficiently close to zero. It would be interesting to investigate whether one can construct higher dimensional counterexamples based on the solutions that were mentioned below (\ref{eqGL}).

For $\varepsilon>0$, let
\begin{equation}\label{eqWeps}
W_\varepsilon(u_1,u_2)=\frac{1}{4}(1-u_1^2-u_2^2)^2+\frac{\varepsilon}{2}u_1^2u_2^2.
\end{equation}
Obviously, this axially symmetric potential has the four points $(\pm 1,0)$, $(0,\pm 1)$ as its only global minimizers. They are nondegenerate but degenerate as $\varepsilon \to 0$. To the best of our knowledge, the investigation of this limiting behavior goes back at least as far as the physics paper \cite{barankov}. Our aforementioned counterexamples to Modica's gradient estimate can be provided at once by combining the next proposition, which is of independent interest, with (\ref{eqIneq}).
\begin{pro}
There exists $\theta_0\in (0,1)$ such that, given $\theta \in (0,\theta_0)$, system (\ref{eqEq}) with $W=W_\varepsilon$, $\varepsilon \in (0,1)$, $n=1$ and  $m=2$  has a $\frac{2\pi}{\theta}$-periodic solution \[U_\varepsilon=(U_{1,\varepsilon},U_{2,\varepsilon}):\mathbb{R}\to \mathbb{R}^2\]
such that
\begin{equation}\label{eqsym1}
 U_{1,\varepsilon}\left(\frac{\pi}{2\theta}-x\right)=U_{2,\varepsilon}(x)>0,\ \ x\in \left(0,\frac{\pi}{2\theta} \right),
\end{equation}
and
\begin{equation}\label{eqsym2}
U_{1,\varepsilon}(-x)= U_{1,\varepsilon}(x),\ \ U_{2,\varepsilon}(-x)=- U_{2,\varepsilon}(x),\ \ x\in \mathbb{R}.
\end{equation}
Moreover, it holds that
\begin{equation}\label{eqconvergence}
\|U_\varepsilon-u_\theta \|_{C^1(\mathbb{R};\mathbb{R}^2)}\to 0\ \ \textrm{as}\ \varepsilon \to 0,
\end{equation}
where $u_\theta$ is as in (\ref{eqper}).
\begin{proof}
Let $\theta>0$. For  $\varepsilon\geq 0$, consider solutions with positive components to the boundary value problem:
\begin{equation}\label{eqBVP}
\left\{
\begin{array}{l}
  u''=\nabla_u W_\varepsilon(u),\ \ x\in \left(0,\frac{\pi}{2\theta} \right); \\
    \\
  u_1'(0)=0,\ u_2(0)=0;\ \ u_1\left(\frac{\pi}{2\theta} \right)=0,\ u_2'\left(\frac{\pi}{2\theta} \right)=0.
\end{array}
\right.
\end{equation}
In view of our previous discussion related to (\ref{eqper}), such a solution exists if $\varepsilon=0$ and is given by
\begin{equation}\label{eqUth}
u_\theta(x)=\left(\sqrt{1-\theta^2}\cos(\theta x),\sqrt{1-\theta^2}\sin(\theta x) \right).
\end{equation}
Moreover, similarly as in \cite[Thm. 4.1]{alamaUniq}, this is the unique solution in this case.
The aforementioned paper dealt with the case of Dirichlet boundary conditions, but the corresponding terms from the integration by parts in the proof still vanish under the above boundary conditions.

By consecutive reflections, a solution to (\ref{eqBVP}) can easily be extended in the whole line as a $\frac{2\pi}{\theta}$-periodic solution with
$u_1$ being even and
$u_2$ odd, as in (\ref{eqsym2}).
The important thing to note here is that, as in \cite{alamaARMA15}, applying the usual maximum priniciple to the scalar function $u_1^2+u_2^2-1$, it follows easily that any such periodic solution should satisfy
\begin{equation}\label{eqcont}
u_1^2(x)+u_2^2(x)<1,\ \ x\in \mathbb{R}.
\end{equation}

We can show the existence of a solution $U_\varepsilon=(U_{1,\varepsilon},U_{2,\varepsilon})$ to (\ref{eqBVP}), satisfying (\ref{eqsym1}), which minimizes the energy
\[
E_\varepsilon(u)=\int_{0}^{\frac{\pi}{2 \theta}}\left[\frac{1}{2}|u'|^2+W_\varepsilon(u) \right]dx
\]
in the space
\[
X=\left\{
\begin{array}{ll}
  u=(u_1,u_2)\in \left[H^1\left(0,\frac{\pi}{2\theta} \right) \right]^2 \ : & u_2(0)=0,\  u_1\left(\frac{\pi}{2\theta} \right)=0, \\
    & \\
    &  \textrm{and}\ u_1\left(\frac{\pi}{2\theta}-x\right)=u_2(x),\  x\in \left(0,\frac{\pi}{2\theta} \right)
\end{array}\right\}.
\]
We claim that there exists $\theta_0>0$ such that  $U_\varepsilon$ is nontrivial if $\theta \in (0,\theta_0)$ and $\varepsilon \in (0,1)$.
Indeed, if not, we can cook up a test function in $X$ which is equal to $(0,1)$ on $\left[1,\frac{\pi}{2\theta}-1 \right]$ to get that
\begin{equation}\label{eqJ0}
J(0,0)=\frac{\pi}{8 \theta}\leq C,
\end{equation}
with $C>0$ independent of both $\theta, \varepsilon \in (0,1)$, which is not possible if $\theta>0$ is sufficiently small and $\varepsilon \in (0,1)$.
We may assume that this solution has positive components because of the property  $E_\varepsilon\left(|u_1|,|u_2|\right)\leq E_\varepsilon(u_1,u_2)$
for any $(u_1,u_2)\in X$, and the strong maximum principle (applied componentwise).

It remains to establish (\ref{eqconvergence}). By virtue of (\ref{eqcont}), using Arzcela-Ascoli's theorem, passing to a subsequence $\varepsilon_j\to 0$, we find that
\[
U_{\varepsilon_j}\to U_0\ \ \textrm{in}\ \ C^1\left(\left[0,\frac{\pi}{2\theta}\right] ; \mathbb{R}^2 \right)\ \ \textrm{as}\ j\to \infty,
\]
where $U_0$ solves (\ref{eqBVP}) with $\varepsilon=0$ and has nonnegative components.
By decreasing $\theta_0$ if necessary,  we can conclude from (\ref{eqJ0}) that $U_0$ is nontrivial. In particular, by the strong maximum priniciple, we deduce that it has positive components.
On the other hand, by our previous discussion in the beginning of the proof, we must have that $U_0=u_\theta$ as in (\ref{eqUth}). Lastly, by the uniqueness
of the limit, we infer that the above convergence holds for $\varepsilon \to 0$, as desired.
\end{proof}
\end{pro}
\begin{rem}
	We observe that system (\ref{eqEq}) with potential as in (\ref{eqWeps}) is symmetric according to the definition of \cite{fazly}. Moreover, it is orientable for positive solutions according to the definition in \cite{fg}.
\end{rem}

\section*{Acknowledgments} We would like to thank Prof. Aftalion for bringing to our attention \cite{barankov} and for useful discussions which motivated the current note. This project has received funding from the European Union's Horizon 2020 research and innovation programme under the Marie Sk\l{}odowska-Curie grant agreement No 609402-2020 researchers: Train to Move (T2M).

\end{document}